\documentclass[a4paper,11pt]{article}
\usepackage[utf8]{inputenc}
\usepackage[T1]{fontenc}

\usepackage{amsthm,amsmath}
\usepackage{tikz}
\usepackage{mathrsfs,amssymb,amsfonts} 
\usepackage{enumitem}
\usepackage{fullpage}
\usepackage{hyperref, enumerate}
\usepackage[babel]{microtype}
\usepackage[english]{babel}
\usepackage[capitalise]{cleveref}

\usepackage{thmtools}
\usepackage{mathtools, comment}
\usepackage{amssymb}
\usepackage[nomath]{lmodern}
\usepackage{graphicx}
\usepackage{pgf,tikz,tkz-graph,subcaption}
\usetikzlibrary{arrows,shapes}
\usetikzlibrary{decorations.pathreplacing}
\usepackage{tkz-berge}
\usepackage{enumitem}
\usepackage[normalem]{ulem}
\usepackage{hyperref}
\hypersetup{colorlinks = true, linkcolor = blue, citecolor = blue, urlcolor = blue}
\newcommand*{\ceilfrac}[2]{\mathopen{}\left\lceil\frac{#1}{#2}\right\rceil\mathclose{}}
\newcommand*{\floorfrac}[2]{\mathopen{}\left\lfloor\frac{#1}{#2}\right\rfloor\mathclose{}}

\allowdisplaybreaks

\usepackage[margin=1in]{geometry}
\newcommand*{\myproofname}{Proof}

\title{Solution of Erd\H{o}s problem $\# 443$}

\author{Stijn Cambie\thanks{Department of Computer Science, KU Leuven Campus Kulak-Kortrijk, 8500 Kortrijk, Belgium. Supported by a postdoctoral fellowship by the Research Foundation Flanders (FWO) with grant number 1225224N. Email: \protect\href{mailto:stijn.cambie@hotmail.com}{\protect\nolinkurl{stijn.cambie@hotmail.com}}}}

\begin{document}
\parindent=0cm
\maketitle

% \begin{abstract}

% \end{abstract}

% \section{Introduction}\label{sec:intro}

In~\cite[p.88]{EG80}, Erd\H{o}s and Graham posed the following problem;

\textit{Let \( m \) and \( n \) be positive integers and consider the two sets 
\[
\left\{ k(m-k) : 1 \leq k \leq \frac{m}{2} \right\} 
\quad \text{and} \quad 
\left\{ \ell(n-\ell) : 1 \leq \ell \leq \frac{n}{2} \right\}.
\]
Can one estimate the number of integers common to both? Is this number unbounded? It should certainly be less than \( (mn)^\varepsilon \) for every \( \varepsilon > 0 \) if \( mn \) is sufficiently large.}

The problem is listed in the online database of Erd\H{o}s problems as problem $\# 443$ (as of the time of publication),~\url{https://www.erdosproblems.com/443}.

We will use the notation \( S(m) = \left\{ k(m-k) : 1 \leq k \leq \frac{m}{2} \right\} \) for such sets, and let $f(n,m) = \# ( S(m) \cap S(n) ).$

Observe that \( f(n,m) \leq f(2n,2m) \) by doubling each of \( n,m,k,\ell \) in an equality $k(m-k)=\ell(n-\ell)$.
Without loss of generality, we will assume $m>n.$

Writing $k(2m-k)=m^2-(m-k)^2=m^2-c^2,$ for $c=m-k$,
we deduce that the number of solutions is bounded by half of the number of divisors of $m^2-n^2$,
since $m^2-c^2=n^2-d^2 \Leftrightarrow m^2-n^2=c^2-d^2=(c-d)(c+d)$ and $c-d$ is thus a divisor bounded by $\sqrt{m^2-n^2}.$

A well-known result~\cite{Ram1915} states that the number of divisors of $x$, $\tau(x)$, is bounded by $2^{(1+o(1)) \frac{\log x}{\log \log x}} = x^{o(1)}$ and this is sharp (for highly composite numbers), up to the $o(1)$ term which is a function tending to zero (as $x \to \infty$).
Since $m^2-n^2<(mn)^2,$ we obtain that $f(m,n) \le f(2m,2n) \le \floorfrac{\tau(m^2-n^2)}{2} = (mn)^{o(1)}$, as expected. 

When $m=n+1$, $f(2m,2n)= \ceilfrac{\tau(2n+1)}{2}-1$ is unbounded, concluding the other question.
For this we note that for every factorization $2n+1=pq$ with $p \ge q\ge 2$ (both odd), $c=\frac{p+q}{2}$ and $d=\frac{p-q}{2}$ are integers for which $c<m$ and $0\le d <n$.

\section*{Remark.}

Norbert Hegyv\'ari~\cite{Heg25} solved the problem independently, $40$ years earlier.
His proof appeared recently, after submission of this note.

% After submission, a solution of the same problem by Hegyv\'ari~\cite{Heg25} has appeared.

\end{document}